\documentclass[a4paper, 10 pt, reqno]{amsart}
\usepackage{amsfonts, amssymb, amsmath, eucal, amsthm}
\usepackage[all]{xy} 

\newtheorem{lemma}{Lemma}
\newtheorem{theorem}[lemma]{Theorem}
\newtheorem{proposition}[lemma]{Proposition}

\newtheorem*{theorem*}{Theorem}

\theoremstyle{definition}
\newtheorem{definition}[lemma]{Definition}
\newtheorem{remark}[lemma]{Remark}
\newtheorem{question}[lemma]{Question}
\newtheorem*{acknowledgments}{Acknowledgments}
\newtheorem{example}[lemma]{Example}

\newcommand{\AdemRuanZhang}{ARZ}
\newcommand{\ChenRuanCohomology}{MR2104605}
\newcommand{\FantechiGottsche}{MR1971293}
\newcommand{\GLSUX}{MR2318652}
\newcommand{\ItoReid}{ MR1463181}
\newcommand{\JKK}{MR2285746}
\newcommand{\AdemLeidaRuan}{MR2359514}
\newcommand{\ChenHu}{MR2242619}

\newcommand{\G}{\mathcal{G}}
\renewcommand{\H}{\mathcal{H}}
\newcommand{\V}{\mathcal{V}}

\title{The age grading and the Chen-Ruan cup product}
\author{Richard Hepworth}
\address{Department of Pure Mathematics\\ University of Sheffield\\  Sheffield\\ S3 7RH}
\email{r.hepworth@sheffield.ac.uk}
\thanks{The author is supported by E.P.S.R.C.~Postdoctoral Research Fellowship EP/D066980.}

\begin{document}

\begin{abstract}We prove that the obstruction bundle used to define the cup-product in Chen-Ruan cohomology is determined by the so-called \emph{age grading} or \emph{degree-shifting numbers}.  Indeed, the obstruction bundle can be directly computed using the age grading.  We obtain a K\"unneth Theorem for Chen-Ruan cohomology as a direct consequence of an elementary property of the age grading, and explain how several other results -- including associativity of the cup-product -- can be proved in a similar way.
\end{abstract}
\maketitle

\section*{Introduction}

In \cite{\ChenRuanCohomology} Chen and Ruan defined the \emph{orbifold} or \emph{Chen-Ruan} cohomology $H^\ast_\mathrm{CR}(X)$ of an almost-complex orbifold $X$.  As a group $H^\ast_\mathrm{CR}(X)$ is simply the cohomology $H^\ast(\Lambda X)$ of the inertia orbifold $\Lambda X$, with degrees shifted by a quantity known as the \emph{age grading} or \emph{degree-shifting number}.  The striking feature of Chen and Ruan's work is that this group is an associative graded ring under the \emph{Chen-Ruan cup product}.  This cup product has proved difficult to compute because its definition involves a certain \emph{obstruction bundle} that Chen and Ruan defined using  orbifold Riemann surfaces.

Since the work of Chen and Ruan several interesting results regarding the obstruction bundle have appeared.  Chen and Hu's de Rham description of the Chen-Ruan cohomology of abelian orbifolds \cite{\ChenHu} involved a computation of the obstruction bundle for abelian orbifolds.  Fantechi and G\"ottsche \cite{\FantechiGottsche} refined Chen and Ruan's construction in the case of global quotient orbifolds.  Jarvis, Kaufmann and Kimura \cite{\JKK} gave an explicit formula for the rational equivariant K-theory class of the Fantechi-G\"ottsche obstruction bundle, so determining it up to isomorphism.

The purpose of this note is to show that the obstruction bundle can be computed directly in terms of the age grading.  As a result we obtain a K\"unneth Theorem for Chen-Ruan cohomology.  We also explain how several known theorems -- including associativity of the Chen-Ruan product -- can be proved as direct consequences of elementary properties of the age grading.

Let us recall from \cite{\ChenRuanCohomology} that the obstruction bundle is a vector-bundle $E\to\Lambda^2X$ over the $2$-sectors of $X$, and is defined using the local description of $\Lambda^2X$ in terms of orbifold-charts on $X$.  We will show that each component $X_{(g_1,g_2)}$ of $\Lambda^2 X$ is naturally equipped with a \emph{twisting group} $\langle g_1, g_2\rangle$  and a fibrewise-linear action of this group on the pullback $\varepsilon^\ast V\to X_{(g_1,g_2)}$ of any vector-bundle $V\to X$.  We obtain the following description of the obstruction bundle.

\begin{theorem}\label{BundleTheorem}
Over a component $X_{(g_1,g_2)}$ of $\Lambda^2 X$ the obstruction bundle $E$ is given by
\[E_{(g_1,g_2)}=(\varepsilon^\ast TX\otimes H^{0,1}_{\bar\partial}(\Sigma))^{\langle g_1,g_2\rangle}\]
where $\Sigma$ is a Riemann surface with action of $\langle g_1,g_2\rangle$ such that $\Sigma/\langle g_1,g_2\rangle$ is an orbifold Riemann-sphere with singular points of order $o(g_1)$, $o(g_2)$ and $o(g_1g_2)$.
\end{theorem}

Chen and Ruan gave a formula for the dimension of the obstruction bundle in terms of the age grading \cite{\ChenRuanCohomology}.  Equipped with Theorem \ref{BundleTheorem}, the study of the obstruction bundle becomes a matter of understanding the assignment
$V\mapsto (V\otimes H^{0,1}_{\bar\partial}(\Sigma))^{\langle g_1,g_2\rangle}$.  By observing that this assignment is determined by the dimension of its values on the irreducible representations of $\langle g_1,g_2\rangle$ we are able to determine the obstruction bundle directly using Chen and Ruan's formula.  Let $\iota_V(g)$ denote the age of $g$ in $V$ and let $V^{g_1,\ldots,g_k}$ denote the elements of $V$ fixed by $g_1,\ldots,g_k$.

\begin{theorem}\label{FunctorTheorem}
Write $V_1,\ldots, V_n$ for the irreducible representations of $\langle g_1,g_2\rangle$ and let $T_i\to X_{(g_1,g_2)}$ be vector bundles for which $\varepsilon^\ast TX= \bigoplus V_i\otimes T_i$ as bundles of $\langle g_1,g_2\rangle$ representations.  Then
\[E_{(g_1,g_2)}=\bigoplus h_i T_i\]
where
\begin{equation}\label{ObstructionDimensionEquation}
h_i=\iota_{V_i}(g_1)+\iota_{V_i}(g_2)-\iota_{V_i}(g_1g_2)+\dim {V_i}^{g_1,g_2}-\dim {V_i}^{g_1g_2}.
\end{equation}
\end{theorem}

\begin{example}\label{Example}
Suppose that $\langle g, h\rangle=\{ \pm 1, \pm g,\pm h,\pm gh\}$ is the quaternion group of order $8$.  Then $E_{(g,h)}$ is the bundle $\mathrm{Hom}_{\langle g, h\rangle}(Q,\varepsilon^\ast TX)$, where $Q$ is the $2$-dimensional irreducible representation of $\langle g,h\rangle$.
\end{example}

Similar methods to those used to prove Theorem~\ref{FunctorTheorem} will be used to recover three existing results.  These are associativity of the Chen-Ruan cup-product \cite{\ChenRuanCohomology}, Chen and Hu's description of the obstruction bundle for abelian orbifolds \cite{\ChenHu}, and Gonz\'alez et al.'s computation of the Chen-Ruan cohomology of cotangent orbifolds \cite{\GLSUX}.  We also obtain the following K\"unneth Theorem:

\begin{theorem}\label{KunnethTheorem}
Let $X$, $Y$ be almost-complex orbifolds with $\mathrm{SL}$ singularities, so that $H^\ast_\mathrm{CR}(X)$ and $H^\ast_\mathrm{CR}(Y)$ are concentrated in integral degrees and we can form the graded ring $H^\ast_\mathrm{CR}(X)\otimes H^\ast_\mathrm{CR}(Y)$.  Then there is a graded ring isomorphism
\[H^\ast_\mathrm{CR}(X)\otimes H^\ast_\mathrm{CR}(Y) \cong H^\ast_\mathrm{CR}(X\times Y).\]
\end{theorem}

\begin{remark}
In de Rham cohomology the cup product can be regarded as the composite 
\[H^\ast(X)\otimes H^\ast(X)\cong H^\ast(X\times X)\xrightarrow{\Delta^\ast}H^\ast(X)\]
of the K\"unneth Isomorphism with the map induced by the diagonal $\Delta\colon X\to X\times X$.  Associativity of the cup-product is then equivalent to $\Delta^\ast(\mathrm{Id}\times\Delta)^\ast=\Delta^\ast(\Delta\times\mathrm{Id})^\ast$, which is just the usual functoriality of induced maps.  Using Theorem~\ref{KunnethTheorem} we can therefore regard Chen and Ruan's definition of the cup-product as a definition of $\Delta^\ast$, and their proof of associativity as a proof that $\Delta^\ast(\mathrm{Id}\times\Delta)^\ast=\Delta^\ast(\Delta\times\mathrm{Id})^\ast$.
\end{remark}

\begin{question}
Is it possible to define the induced map $f^\ast\colon H^\ast_\mathrm{CR}(Y)\to H^\ast_\mathrm{CR}(X)$ associated to a general map of orbifolds $f\colon X\to Y$?  Does functoriality $g^\ast f^\ast=(fg)^\ast$ hold?
\end{question}

\begin{remark}
One might wish to take Theorems \ref{BundleTheorem} and \ref{FunctorTheorem} as the \emph{definition} of the obstruction bundle, so removing the need for orbifold Riemann surfaces.  This is not possible since the proof of Theorem \ref{FunctorTheorem} requires an application of Chen and Ruan's formula \cite{\ChenRuanCohomology} for the dimension of the obstruction bundle in order to show that the right-hand-side of \eqref{ObstructionDimensionEquation} is non-negative.  Nevertheless, after this single appeal to the theory of orbifold Riemann surfaces, whose conclusion is simply an inequality regarding the age grading, one can regard the results here as an elementary way to define the Chen-Ruan cup product and to prove its associativity.  A positive answer to either part of the following question would remove the need for orbifold Riemann surfaces entirely.
\end{remark}

\begin{question}
Is there an elementary proof of the inequality
\[\iota_V(g_1)+\iota_V(g_2)-\iota_V(g_1g_2)+\dim V^{g_1,g_2}-\dim V^{g_1g_2}\geqslant 0?\]
Is there an elementary description of the assignment
\[V\mapsto (V\otimes H^{0,1}_{\bar\partial}(\Sigma))^{\langle g_1,g_2\rangle}?\]
\end{question}

As mentioned earlier Jarvis, Kaufmann and Kimura \cite{\JKK} have determined the obstruction bundle of a global quotient by giving an explicit description of the rational $K$-theory class of the Fantechi-G\"ottsche obstruction bundle.  Their methods, combined with Theorem \ref{BundleTheorem}, could be used to give an analogous result for general orbifolds.  The main step -- corresponding to \cite[Lemma 8.5]{\JKK} -- in the proof of such a result  would rely on the same key observation that is used in the proof of Theorem \ref{FunctorTheorem}: that Chen and Ruan's formula for the dimension of the obstruction bundle in fact determines the obstruction bundle entirely.  The same observation was also used in the proof of \cite[Proposition 3.4]{\ChenHu}.

Here is an outline of the paper.  In Section~\ref{TautologicalSection} we recall the definition of the $k$-sectors $\Lambda^kX$ of $X$.  We characterize $\Lambda^k X$ in a way that allows us to introduce, for each component $X_{(g_1,\ldots,g_k)}$, a \emph{twisting group} $\langle g_1,\ldots,g_k\rangle$ and its \emph{twisting action} on the pullback to $X_{(g_1,\ldots,g_k)}$ of any vector-bundle $V\to X$.  We then prove Theorem~\ref{BundleTheorem}.  In Section~\ref{AgeObstructionSection} we recall the age grading and list some of its properties.  We then prove Theorem~\ref{FunctorTheorem}.  Section~\ref{ExampleSection} shows how to perform the calculation stated in Example~\ref{Example}. In Section~\ref{ApplicationSection} we prove Theorem~\ref{KunnethTheorem} and outline how similar methods may be used in the proof of the results of Chen and Ruan \cite{\ChenRuanCohomology}, Chen-Hu \cite{\ChenHu}, and Gonz\'alez et.~al \cite{\GLSUX} mentioned earlier.

\begin{acknowledgments}
The author is supported by an E.P.S.R.C.~Postdoctoral Research Fellowship, grant number EP/D066980.
\end{acknowledgments}

\section{Twisted sectors and the {twisting} group}\label{TautologicalSection}

To an orbifold $X$ one can associate the \emph{$k$-sectors} or \emph{twisted $k$-sectors} of $X$, which are orbifolds $\Lambda^k X$ for $k\geqslant 0$.  
The $0$-sectors $\Lambda^0 X$ is just $X$ itself.  The $1$-sectors $\Lambda X:=\Lambda^1 X$ is called the \emph{inertia orbifold}.  There is an \emph{evaluation map}
\[\varepsilon\colon\Lambda^k X\to X,\]
and more general evaluation maps $\varepsilon_{I_1\ldots I_j}\colon\Lambda^kX\to\Lambda^j X$ for any sequence $I_1,\ldots,I_j$ of ordered tuples in $\{1,\ldots,k\}$.  See \cite[\S 2]{\AdemRuanZhang}, or \cite[\S 4.1]{\AdemLeidaRuan}.

The following proposition gives a new characterization of the twisted sectors, in terms of which we can immediately define the twisting group and the twisting action.  The section ends with the proof of Theorem \ref{BundleTheorem}.  In what follows we will not distinguish between the orbifold $X$ and the groupoid $\G$ representing it; our constructions are Morita-invariant.

\begin{proposition}\label{TwistedSectorsProposition}
Let $\G$ be a proper \'etale Lie groupoid and $\H$ a Lie groupoid.
\begin{enumerate}
\item\label{PartOne} Morphisms $\H\to\Lambda^k\G$ correspond precisely to diagrams of the form 
\[\xymatrix{
\H\ar@/^1.5pc/[rrr]_{\quad}^{f}_{}="1"\ar@/_1.5pc/[rrr]_{f}^{}="2" &&& \G.\ar@{=>}"1";"2"^{\phi_1,\ldots,\phi_k}
}\]
The morphism corresponding to the diagram above will be written as $(f,\phi_i)$.
\item\label{PartTwo} $2$-morphisms 
\[\xymatrix{
\H\ar@/^1.5pc/[rr]^{(f,\phi_i)}_{}="1"\ar@/_1.5pc/[rr]_{(g,\gamma_i)}^{}="2" &&\Lambda^k\G\ar@{=>}^{\psi}"1";"2"
}\]
correspond precisely to $2$-morphisms $\psi\colon f\Rightarrow g$ for which $\psi\phi_i=\gamma_i\psi$.
\item\label{PartThree}
$\mathrm{Id}\colon\Lambda^k\G\to\Lambda^k\G$ corresponds to
\[\xymatrix{
\Lambda^k\G\ar@/^1.5pc/[rrr]^{\varepsilon}_{}="1"\ar@/_1.5pc/[rrr]_{\varepsilon}^{}="2"&&&\G\ar@{=>}^{E_1,\ldots,E_k}"1";"2"
}\]
where $\varepsilon$ is the usual evaluation map and the $E_i$ are canonically-determined $2$-automorphisms of $\varepsilon$.
\item\label{PartFour}
$\varepsilon_{I_1\ldots I_j}\colon\Lambda^k\G\to\Lambda^j\G$ corresponds to
\[\xymatrix{
\Lambda^k\mathcal{G}\ar@/^1.5pc/[rrr]^{\varepsilon}_{}="1"\ar@/_1.5pc/[rrr]_{\varepsilon}^{}="2"&&&\mathcal{G}\ar@{=>}^{E_{I_1},\ldots,E_{I_j}}"1";"2"
}\]
where $E_{I_l}=E_{l_1}\circ\cdots\circ E_{l_m}$ for $I_l=(l_1,\ldots,l_m)$.
\end{enumerate}
\end{proposition}

\begin{definition}We will denote the components of $\Lambda^kX$ by $X_{(g_1,\ldots,g_k)}$, $X_{(h_1,\ldots,h_k)}$, etcetera; the $g_i$ and the $h_i$ are simply labels for the components.
\begin{enumerate}
\item Let $X_{(g_1,\ldots,g_k)}$ be a component of $\Lambda^k X$.  The \emph{twisting group} $\langle g_1,\ldots,g_k\rangle$ of $X_{(g_1,\ldots,g_k)}$ is the group on generators $g_1,\ldots,g_k$, isomorphic under $g_i\mapsto E_i$ to $\langle E_1,\ldots,E_k\rangle\subset\mathrm{Aut}(\varepsilon|_{X_{(g_1,\ldots,g_k)}}\colon X_{(g_1,\ldots,g_k)}\to X)$.
\item Let $V\to X$ be a vector bundle.  The \emph{tautological action} of $\langle g_1,\ldots, g_k\rangle$ on the bundle $\varepsilon^\ast V\to X_{(g_1,\ldots,g_k)}$ is given by the fibrewise-linear automorphisms $g_i\colon \varepsilon^\ast V\to\varepsilon^\ast V$ induced by the $2$-automorphisms $E_i$ of $\varepsilon$.
\end{enumerate}
\end{definition}

Finiteness of the twisting group is guaranteed by Lemma \ref{DMAutLemma} below.  For the second part of the definition recall that, given $f\colon\mathcal{A}\to\mathcal{C}$ and $g\colon\mathcal{B}\to\mathcal{C}$, $2$-automorphisms of $f$ induce automorphisms of $\mathcal{A}\times_{f,g}\mathcal{B}$, and that when $g\colon\mathcal{B}\to\mathcal{C}$ is a vector-bundle these automorphisms are fibrewise-linear maps.  The twisting group and twisting action also satisfy certain naturality properties with respect to maps of $X$, maps of $V$, and evaluation maps, all of which follow from Proposition \ref{TwistedSectorsProposition}.  

\begin{remark}
Proposition \ref{TwistedSectorsProposition} suggests a definition of the twisted sectors $\Lambda^k\mathfrak{X}$ of a differentiable Deligne-Mumford stack $\mathfrak{X}$, where one defines morphisms $U\to\Lambda^k\mathfrak{X}$ using an analogue of Proposition \ref{TwistedSectorsProposition}, part \ref{PartOne}, and one defines $2$-morphisms and evaluation maps using analogues of parts \ref{PartTwo}, \ref{PartThree}, and \ref{PartFour}.  This is indeed possible, and one finds that if $\G$ is a groupoid representing $\mathfrak{X}$, then $\Lambda^k\G$ is precisely the groupoid representing $\mathfrak{X}$ that one obtains from $\G$ by pulling back under $\varepsilon\colon\Lambda^k\mathfrak{X}\to \mathfrak{X}$.  When $k=1$ the resulting inertia stack $\Lambda\mathfrak{X}$ is equivalent to, but leaner than, the more usual $\Lambda\mathfrak{X}:=\mathfrak{X}\times_{\mathfrak{X}\times\mathfrak{X}}\mathfrak{X}$.
\end{remark}

\begin{lemma}\label{DMAutLemma}
Suppose given a diagram of the form 
\[\xymatrix{
a\ar[r] & \mathcal{A}\ar@/^1.5pc/[rr]^{m}_{}="1"\ar@/_1.5pc/[rr]_{m}^{}="2" &&\mathcal{B}\ar@{=>}^{\phi}"1";"2"
}\]
where $a$ is a point in a \emph{connected} Lie groupoid $\mathcal{A}$, and $\mathcal{B}$ is a proper \'etale Lie groupoid.  Then $\phi$ is trivial if and only if $\phi|_a\colon m|_a\Rightarrow m|_a$ is trivial.
\end{lemma}
\begin{proof}
This is immediate in the case $\mathcal{B}=U\rtimes G$ for $G$ a finite group acting on $U$, and follows in general because $\mathcal{B}$ locally has this form.
\end{proof}

\begin{proof}[Proof of Proposition \ref{TwistedSectorsProposition}]
Recall from \cite{\AdemRuanZhang} that $\Lambda^k\G$ is the groupoid with objects and arrows
\begin{eqnarray*}
\Lambda^k\G_0&=&\{ (a_1,\ldots,a_k)\in\G_1^k\mid s(a_i),\ t(a_j)\ \mathrm{all\ coincide}\},\\
\Lambda^k\G_1&=&\{ (u,a_1,\ldots,a_k)\in\G_1^{k+1}\mid s(a_i),\ t(a_j),\ s(u)\ \mathrm{all\ coincide}\},
\end{eqnarray*}
with source and target
\begin{align*}
s(u,a_1,\ldots,a_k)&=(a_1,\ldots,a_k),\\
t(u,a_1,\ldots,a_k)&=(ua_1u^{-1},\ldots,ua_ku^{-1}),
\end{align*}
and structure-maps
\begin{gather*}
e(a_1,\ldots,a_k)=(e(\alpha),a_1,\ldots,a_k),\ \ \ \alpha=s(a_i)=t(e_j),\\
i(u,a_1,\ldots,a_k)=(u^{-1},ua_1u^{-1},\ldots,ua_ku^{-1}),\\
m((v,ua_1u^{-1},\ldots,ua_ku^{-1}),(u,a_1,\ldots,a_k))=(vu,a_1,\ldots,a_k).
\end{gather*}

Suppose given $f\colon\H\to\G$ and $\phi_i\colon f\Rightarrow f$, $i=1,\ldots,k$.  We obtain
\begin{align*}
F_0\colon\H_0&\to\G_1^k,\\
F_1\colon\H_1&\to\G_1^{k+1},
\end{align*}
given by
\begin{align*}
F_0(h)&=(\phi_1(h),\ldots,\phi_k(h)),\\
F_1(H)&=(f_1(H),\phi_1(sH),\ldots,\phi_k(sH)).
\end{align*}
We claim that these maps together give a groupoid morphism $\H\to\Lambda^k\G$.  That $F_0$, $F_1$ are maps into $\Lambda^k\G_0$, $\Lambda^k\G_1$ respectively follows from $s\circ\phi_i=f_0=t\circ\phi_i$.  That $F_0$, $F_1$ respect the source and target maps follows from the same fact for $f$, together with the naturality of the $\phi_i$.  That $F_1$ commutes with composition follows from the same fact for $f_1$ together with the naturality of the $\phi_i$.  Thus we have constructed a groupoid morphism $F\colon\H\to\Lambda^k\G$ as required.  This reasoning can be reversed to produce $f\colon\H\to\G$, $\phi_i\colon f\Rightarrow f$ from a groupoid-morphism $F\colon\H\to\Lambda^k\G$, and part \ref{PartOne} is proved.  The proof of part \ref{PartTwo} is similar.

The evaluation map $\varepsilon\colon\Lambda^k\G\to\G$ is the map obtained by sending the $a_i$ to their common source and target.  From the proof of part \ref{PartOne} it is therefore immediate that $\Lambda^k\G\to\Lambda^k\G$ corresponds to 
\[\xymatrix{
\Lambda^k\G\ar@/^1.5pc/[rrr]^{\varepsilon}_{}="1"\ar@/_1.5pc/[rrr]_{\varepsilon}^{}="2"&&&\G\ar@{=>}^{E_1,\ldots,E_k}"1";"2"
}\]
where $E_i\colon\Lambda^k\G_0\to\G_1$ is $(a_1,\ldots,a_k)\mapsto a_i$.  The diagram
\[\xymatrix{
\Lambda^k\mathcal{G}\ar@/^1.5pc/[rrr]^{\varepsilon}_{}="1"\ar@/_1.5pc/[rrr]_{\varepsilon}^{}="2"&&&\mathcal{G}\ar@{=>}^{E_{I_1},\ldots,E_{I_j}}"1";"2"
}\]
induces $\Lambda^k\G\to\Lambda^j\G$, $(a_1,\ldots,a_k)\mapsto(a_{I_1},\ldots,a_{I_j})$ on objects, and similarly for arrows, where $a_{I_l}=a_{l_1}\cdots a_{l_m}$, $I_l=(l_1,\ldots,l_m)$; this is just the evaluation map $\varepsilon_{I_1\cdots I_j}$.  Parts \ref{PartThree} and \ref{PartFour} are proved.
\end{proof}

\begin{proof}[Proof of Theorem \ref{BundleTheorem}]
Consider an orbifold groupoid of the form $U\rtimes G$, where $G$ is a finite group acting on a manifold $U$.  Then $\Lambda^k(U\rtimes G)=\bigsqcup U^{h_1,\ldots,h_k}\rtimes G$, where the union is taken over all $k$-tuples in $G$.  The evaluation map $\varepsilon\colon\Lambda^k(U\rtimes G)\to U\rtimes G$ is induced by the componentwise inclusion $\bigsqcup U^{h_1,\ldots,h_k}\to U$, and $E_i\colon \bigsqcup U^{h_1,\ldots,h_k}\to U\times G$ sends $u\in U^{h_1,\ldots,h_k}$ to $(u,h_i)$.  A component of $\Lambda^k(U\rtimes G)$ is then a component of 
\begin{eqnarray*}U\rtimes G_{(g_1,\ldots,g_k)}&=&\bigsqcup U^{h_1,\ldots,h_k}\rtimes G\\ &\simeq& U^{h_1,\ldots,h_k}\rtimes C_G(h_1,\ldots,h_k)\end{eqnarray*}
where $(g_1,\ldots,g_k)$ is a diagonal conjugacy class in $G$; on the first line the union runs over $(h_1,\ldots,h_k)\in(g_1,\ldots,g_k)$ and on the second line a single choice of such an $(h_1,\ldots,h_k)$ has been made.  Then the twisting group $\langle g_1,\ldots,g_k\rangle$ is isomorphic under $g_i\mapsto h_i$ to $\langle h_1,\ldots,h_k\rangle\subset G$.

Now let $V\to U\rtimes G$ be a vector-bundle, which is to say that $V$ is a $G$-equivariant bundle over $U$.  Then at any point $u\in U^{h_1,\ldots,h_k}$, $\langle h_1,\ldots,h_k\rangle$ acts on the fibre of $V$ at $U$.  Under $g_i\mapsto h_i$ this is precisely the twisting action of $\langle g_1,\ldots,g_k\rangle$ on the fibres of $\varepsilon^\ast V\to U\rtimes G_{(g_1,\ldots,g_k)}$.

Recall from \cite[\S 4.3]{\AdemLeidaRuan} the construction of the obstruction-bundle $E_{(g_1,g_2)}\to\G_{(g_1,g_2)}$ over a component of $\Lambda^2\G$.  Let $(y,h_1,h_2)$ be a point of $\G_{(g_1,g_2)}$.  Take an orbifold-chart $U_y\rtimes G_y$ around $y$ in $\G$, so that $\bigsqcup_{(h'_1,h'_2)\sim(h_1,h_2)}U_y^{h'_1,h'_2}\rtimes G_y\simeq U_y^{h_1,h_2}\rtimes C_G(h_1,h_2)$ is an orbifold-chart around $(y,h_1,h_2)$ in $\G_{(g_1,g_2)}$.  Let $N_y=\langle h_1,h_2\rangle\subset G_y$.  Consider the pullback tangent-bundle $\varepsilon^\ast T\G\to\G_{(g_1,g_2)}$ and the vector-space $H^{0,1}_{\bar\partial}(\Sigma_y)$, where $\Sigma_y$ is the Riemann surface with $N_y$-action such that $\Sigma_y/N_y$ is the orbifold Riemann sphere with marked points of order $o(h_1)$, $o(h_2)$, $o(h_2^{-1}h_1^{-1})$, respectively.  Then $N_y$ acts on both $H^{0,1}_{\bar\partial}(\Sigma_y)$ and $\varepsilon^\ast T\G$, and over the chosen chart for $\G_{(g_1,g_2)}$ the obstruction bundle is defined to be 
\[(H^{0,1}_{\bar\partial}(\Sigma_y)\otimes\varepsilon^\ast T\G)^{N_y}.\]
Allowing $y$ to vary, one obtains $E_{(g_1,g_2)}\to\G_{(g_1,g_2)}$.  

Now let $\Sigma$ be the Riemann surface with $\langle g_1,g_2\rangle$-action required for the theorem.  Note, using the first paragraph, that $g_i\mapsto h_i$ identifies the twisting group $\langle g_1,g_2\rangle$ with $N_y$, so that we may take $\Sigma_y=\Sigma$.  Note also that under $N_y\cong\langle g_1,g_2\rangle$, the action of $N_y$ on $\varepsilon^\ast T\G$ is just the tautological action.  Thus, over the chosen orbifold-chart for $(y,h_1,h_2)$, the last paragraph states that the obstruction-bundle is
\[(H^{0,1}_{\bar\partial}(\Sigma)\otimes\varepsilon^\ast T\G)^{\langle g_1,g_2\rangle}.\]
By allowing $y$ to vary, the theorem is proved.
\end{proof}

\section{The age grading.}\label{AgeObstructionSection}

In this section we recall the \emph{age-grading} or \emph{degree-shifting numbers} and we list some properties.  We then discuss  \emph{additive functors} and \emph{additive functions} before proving Theorem~\ref{FunctorTheorem}.

\begin{definition}\label{AgeDefinition}
Let $G$ be a finite group, $g$ an element of $G$, and $V$ a complex representation of $G$.  The \emph{age of g}, denoted by $\iota_V(g)$, is
\[\iota_V(g)=\sum\lambda_i,\]
where $g$ has the form
\[\left(\begin{array}{ccc}e^{2\pi i\lambda_1}&&\\ &\ddots&\\ && e^{2\pi i\lambda_n}\end{array}\right),\qquad 0\leqslant\lambda_i<1\]
with respect to an appropriate basis of $V$.  The age grading appears in \cite{\ItoReid} for $G\subset\mathrm{SL}(n,\mathbb{C})$ and $V=\mathbb{C}^n$. Chen and Ruan \cite{\ChenRuanCohomology} gave the slightly more general definition above under the name \emph{degree-shifting number}.
\end{definition}

\begin{lemma}\label{AgeLemma}\hfill
\begin{enumerate}
\item $\iota_V(g)$ depends only on the isomorphism class of $V$ and the conjugacy class of $g$.
\item $\iota_{V\oplus W}(g)=\iota_V(g)+\iota_W(g)$.
\item $\exp(2\pi i \iota_V(g))=\det(g\colon V\to V)$.
\item $\iota_V(g_1)+\iota_V(g_2)-\iota_V(g_1g_2)+\dim V^{g_1,g_2}-\dim V^{g_1g_2}\geqslant 0$.
\end{enumerate}
\end{lemma}
\begin{proof}
It is trivial to verify the first three properties.  The last is due to Chen and Ruan: the statement in \cite[Theorem 4.1.5]{\ChenRuanCohomology}, that the orbifold cup-product preserves the grading of $H^\ast_\mathrm{CR}(X)$, when applied to the orbifold $X=V/G$, is precisely the statement that the left-hand-side is the dimension of the obstruction bundle, and so is a non-negative integer.  This can also be deduced from the Riemann-Roch formula as in \cite[\S 1]{\FantechiGottsche}, and from the Eichler trace formula as explained in \cite[\S 8]{\JKK}.
\end{proof}

To prove Theorem~\ref{FunctorTheorem} we must consider the assignment $V\mapsto (V\otimes H^{0,1}_{\bar\partial}(\Sigma))^{\langle g_1,g_2\rangle}$, which we regard as a functor from representations of $\langle g_1,g_2\rangle$ to vector spaces.  More generally we shall consider \emph{additive} functors $\V_G\to\V$.  Here $\V$ is the category of finite-dimensional complex vector-spaces, $\V_G$ is the category of finite-dimensional complex representations of a finite group $G$, and additive means that the functor preserves direct sums.  We shall also consider \emph{additive} functions $|\V_G|\to\mathbb{N}\cup\{0\}$, where $|\V_G|$ denotes the isomorphism classes in $\V_G$ and {additive} means that the function sends direct sums to sums.  Each additive functor $H$ yields an additive function $\dim_H$ by taking the dimension, and any additive function $f$ arises in this way by setting $H_f(-)=\bigoplus_i f(V_i)\hom_G(V_i,-)$, where the $V_i$ are the distinct irreducible representations of $G$.  By basic representation theory this establishes a $1-1$ correspondence between the additive functions and natural isomorphism classes of additive functors.

Our decision to consider additive functors, rather than the representations that afford them, will be justified in the next section, where we will have to consider functors such as $V\mapsto V^g$.   These functors are easy to write down, but the representations that afford them are not.

\begin{proof}[Proof of Theorem \ref{FunctorTheorem}]
Chen and Ruan's result \cite{\ChenRuanCohomology} that the orbifold cup-product preserves the grading of $H^\ast_{\mathrm{CR}}(X)$, when applied to the orbifold $V/\langle g_1,g_2\rangle$, states precisely that $(V\otimes H^{0,1}_{\bar\partial}(\Sigma))^{\langle g_1,g_2\rangle}$ has dimension
\begin{equation}\label{BigEquation}\iota_V(g_1)+\iota_V(g_2)-\iota_V(g_1g_2)+\dim V^{g_1,g_2}-\dim V^{g_1g_2}.\end{equation}
The assignment $V\mapsto (V\otimes H^{0,1}_{\bar\partial}(\Sigma))^{\langle g_1,g_2\rangle}$ is an additive functor $\V_{\langle g_1,g_2 \rangle}\to\V$ and the corresponding additive function sends $[V]$ to the expression \eqref{BigEquation}.  But the assignment $V\mapsto\bigoplus h_i\mathrm{Hom}_{\langle g_1,g_2\rangle}(V_i,V)$ is a second additive functor that corresponds to this additive function, by the second part of Lemma~\ref{AgeLemma}.  Consequently the two functors are naturally isomorphic, and so
\[E_{(g_1,g_2)}=(H^{0,1}_{\bar\partial}(\Sigma)\otimes \varepsilon^\ast TX)^{\langle g_1,g_2\rangle}\cong\bigoplus h_i\mathrm{Hom}(V_i,\varepsilon^\ast TX)=\bigoplus h_i T_i.\]
This completes the proof.
\end{proof}

\section{An example}\label{ExampleSection}

Suppose that we wish to compute the obstruction bundle over a $2$-sector $X_{(g,h)}$ of an orbifold $X$, and that the twisting group $\langle g,h\rangle$ is the quaternion group of order $8$.  Thus
\[\langle g,h\rangle = \left\{ \pm 1, \pm g, \pm h,\pm gh\right\}\]
where $-1$ is central, $-k$ denotes $-1\cdot k$, and $g^2=h^2=-1$ and $gh=-hg$.

To begin we must find the irreducible representations of $\langle g,h\rangle$.  These are 
\[ 1,\ G,\ H,\ GH,\ Q,\]
where $1$ is the trivial representation, 
$G$ is the linear representation on which $g=-1$ and $h=1$,
$H$ is the linear representation on which $g=1$ and $h=-1$,
$GH$ is the linear representation on which $g=h=-1$, and  $Q$ is the $2$-dimensional representation on which
\[g=\left(\begin{array}{cc} i & 0 \\ 0 & -i \end{array}\right)\quad \mathrm{and}\quad h=\left(\begin{array}{cc} 0 & 1 \\ -1 & 0 \end{array}\right).\]

Now we must compute the quantities
\[h_V=\iota_{V}(g)+\iota_{V}(h)-\iota_{V}(gh)+\dim {V}^{g,h}-\dim {V}^{gh}\]
for each irreducible represenation $V$.  We find that
\begin{eqnarray*}
h_1&=&0+0-0+0-0\\
h_G&=&{\textstyle\frac{1}{2}}+0-{\textstyle\frac{1}{2}}+0-0\\
h_H&=&0+{\textstyle\frac{1}{2}}-{\textstyle\frac{1}{2}}+0-0\\
h_{GH}&=&{\textstyle\frac{1}{2}}+{\textstyle\frac{1}{2}}-0+0-1\\
h_{Q}&=&1+1-1+0-0
\end{eqnarray*}
so that the $h_V$ all vanish except for $h_Q$, which is equal to $1$.  Now we can apply Theorem~\ref{FunctorTheorem} and compute $E_{(g,h)}$.  Note that in the theorem $T_i=\mathrm{Hom}_{\langle g_1,g_2\rangle}(V_i,\varepsilon^\ast TX)$, and so immediately we obtain
\[E_{(g,h)}=\mathrm{Hom}_{\langle g,h\rangle}(Q,\varepsilon^\ast TX)\]
as claimed in Example~\ref{Example}.

\section{Applications}\label{ApplicationSection}

In this section we prove Theorem~\ref{KunnethTheorem} using an elementary property of the age grading.  We then explain how a similar method can be used in proofs of Chen and Ruan's result on the associativity of the cup-product  \cite{\ChenRuanCohomology}, Chen and Hu's computation of the obstruction bundle of abelian orbifolds  \cite{\ChenHu}, and Gonz\'alez et al.'s computation of the Chen-Ruan cohomology of cotangent orbifolds \cite{\GLSUX}.

\begin{proof}[Proof of Theorem \ref{KunnethTheorem}]
To begin with we note that, given complex representations $V$ of $G$ and $W$ of $H$, and elements $g\in G$, $h\in H$, we have
\begin{equation}\label{AgeProductEquation}\iota_{V\oplus W}(g\times h)=\iota_V(g)+\iota_W(h).\end{equation}

Proposition \ref{TwistedSectorsProposition} gives an isomorphism $\Lambda^k(X\times Y)\cong\Lambda^k X\times\Lambda^k Y$ under which $\varepsilon\colon\Lambda^k(X\times Y)\to X\times Y$ and its $2$-automorphisms $E_i$ correspond to $\varepsilon\times\varepsilon$ and $E_i\times E_i$ respectively.  We can therefore write $(X\times Y)_{(g_1\times h_1,\ldots,g_k\times h_k)}$ for the component $X_{(g_1,\ldots,g_k)}\times Y_{(h_1,\ldots,h_k)}$, where $\langle g_1\times h_1,\ldots,g_k\times h_k\rangle$ is indeed identified with the subgroup of $\langle g_1,\ldots,g_k\rangle\times\langle h_1,\ldots, h_k\rangle$ generated by the $g_i\times h_i$.   Furthermore, $\varepsilon^\ast T(X\times Y)\to (X\times Y)_{(g_1\times h_1,\ldots,g_k\times h_k)}$ is equivariantly identified with $\varepsilon^\ast TX\oplus\varepsilon^\ast TY$.

Using the last paragraph we have the usual K\"unneth Isomorphism $H^\ast(\Lambda X)\otimes H^\ast(\Lambda Y)\cong H^\ast(\Lambda(X\times Y))$, and by \eqref{AgeProductEquation} and the last paragraph the grading shifts respect this isomorphism, so that we have an isomorphism of graded vector-spaces $H^\ast_\mathrm{CR}(X)\otimes H^\ast_\mathrm{CR}(Y)\cong H^\ast_\mathrm{CR}(X\times Y)$.  We must prove that this is  a ring-homomorphism.

Let $\tau_1,\tau_2\in H^\ast_\mathrm{CR}(X)$, $\sigma_1,\sigma_2\in H^\ast_\mathrm{CR}(Y)$.  Then
\begin{eqnarray*}(\tau_1\cup_\mathrm{CR}\tau_2)\times(\sigma_1\cup_\mathrm{CR}\sigma_2)
&=&{\varepsilon_{12}}_\ast(\varepsilon_1^\ast\tau_1\cup\varepsilon_2^\ast\tau_2\cup e(E))\times {\varepsilon_{12}}_\ast(\varepsilon_1^\ast\sigma_1\cup\varepsilon_2^\ast\sigma_2\cup e(E))\\
&=&{\varepsilon_{12}}_\ast\left( (\varepsilon_1^\ast\tau_1\cup\varepsilon_2^\ast\tau_2\cup e(E))\times(\varepsilon_1^\ast\sigma_1\cup\varepsilon_2^\ast\sigma_2\cup e(E))\right)\\
&=&(-1)^d{\varepsilon_{12}}_\ast(\varepsilon_1^\ast(\tau_1\times\sigma_1)\cup\varepsilon_2^\ast(\tau_2\times\sigma_2)\cup(e(E)\times e(E)))\\
&=&(-1)^d(\tau_1\times\sigma_1)\cup_\mathrm{CR}(\tau_2\times\sigma_2),
\end{eqnarray*}
where $d={\mathrm{deg}(\tau_2)\cdot \mathrm{deg}(\sigma_1)}$, as required.  The third line holds because the Euler classes have even degrees and -- since $X$ and $Y$ have $\mathrm{SL}$ singularities -- the honest degrees (as elements of $H^\ast(\Lambda X)$, $H^\ast(\Lambda Y)$) of $\tau_2$ and $\sigma_1$ agree with their shifted degrees modulo 2.  The last line relies on an isomorphism
\begin{equation}\label{KunnethEquation}\pi_1^\ast E\oplus\pi_2^\ast E\cong E\end{equation}
of obstruction bundles over $(X\times Y)_{(g_1\times h_1,g_2\times h_2)}\cong X_{(g_1,g_2)}\times Y_{(h_1,h_2)}$.  Let us write $H_{(g_1,g_2)}$ for the functor $V\mapsto(V\otimes H^{0,1}_{\bar\partial}(\Sigma))^{\langle g_1,g_2 \rangle}$.  Then the required isomorphism \eqref{KunnethEquation} will follow, using Theorem \ref{BundleTheorem}, Theorem \ref{FunctorTheorem} and the comments at the start of the proof, from a natural isomorphism
\[H_{(g_1\times h_1,g_2\times h_2)}(V\oplus W)\cong H_{(g_1,g_2)}(V)\oplus H_{(h_1,h_2)}(W)\]
of functors $\mathcal{V}_G\times\mathcal{V}_H\to\mathcal{V}$.  Since the functors are additive, this will follow from natural isomorphisms
\begin{align*}
H_{(g_1\times h_1,g_2\times h_2)}(V\oplus 0)&\cong H_{(g_1,g_2)}(V),\\
H_{(g_1\times h_1,g_2\times h_2)}(0\oplus W)&\cong H_{(h_1,h_2)}(W).
\end{align*}
But by \eqref{AgeProductEquation} the dimensions of the left hand sides are equal to the dimensions of the right hand sides.  The isomorphisms now follow from the discussion in Section \ref{AgeObstructionSection}.
\end{proof}

Now we shall sketch the proofs of three other well-known results, explaining in each case how one can simplify the proof using the techniques presented in this paper.

\begin{theorem}[Chen-Ruan \cite{\ChenRuanCohomology}]
The Chen-Ruan cup product is associative.
\end{theorem}

\begin{proof}[Sketch Proof.]
Chen and Ruan's proof of this result combines a cohomological argument with \cite[Lemma 4.3.2]{\ChenRuanCohomology}, the essential consequence of which is that there is an isomorphism of vector-bundles:
\begin{equation}\label{AssociativityIsomorphism}\epsilon_{12,3}^\ast E\oplus \epsilon_{1,2}^\ast E\oplus \mathrm{Exc}_{12} \cong \epsilon_{1,23}^\ast E\oplus \epsilon_{2,3}^\ast E\oplus \mathrm{Exc}_{23}.\end{equation}
Here $\epsilon_{12,3}$, $\epsilon_{1,2}$, $\epsilon_{1,23}$, $\epsilon_{2,3}$ are evaluation-maps $\Lambda^3X\to\Lambda^2X$, and $\mathrm{Exc}_{12}$, $\mathrm{Exc}_{23}$ are the `excess bundles' obtained by applying the functors  $V\mapsto V^{g_1g_2}/(V^{g_1,g_2}+V^{g_1g_2,g_3})$, $V\mapsto V^{g_2g_3}/(V^{g_2,g_3}+V^{g_1,g_2g_3})$ to $\epsilon^\ast TX$.

Equation \eqref{AssociativityIsomorphism} was proved in \cite{\ChenRuanCohomology} by manipulating orbifold Riemann-surfaces.  Here we shall show how it follows by considering the age grading.  Over a fixed component $X_{(g_1,g_2,g_3)}$ of $\Lambda^3 X$ each side of \eqref{AssociativityIsomorphism} is obtained by applying an additive functor to $\varepsilon^\ast TX$.  The two corresponding additive functions send a representation $V$ to the quantities
\begin{multline*}
\iota_V(g_1g_2)+\iota_V(g_3)-\iota_V(g_1g_2g_3)-\dim V^{g_1g_2g_3}+\dim V^{g_1g_2,g_3}\\
+\iota_V(g_1)+\iota_V(g_2)-\iota_V(g_1g_2)-\dim V^{g_1g_2}+\dim V^{g_1,g_2}\\
+\dim V^{g_1g_2}-\dim V^{g_1,g_2}-\dim V^{g_1g_2,g_3}+\dim V^{g_1,g_2,g_3}
\end{multline*}
and
\begin{multline*}
\iota_V(g_1)+\iota_V(g_2g_3)-\iota_V(g_1g_2g_3)-\dim V^{g_1g_2g_3}+\dim V^{g_1,g_2g_3}\\
+\iota_V(g_2)+\iota_V(g_3)-\iota_V(g_2g_3)-\dim V^{g_2g_3}+\dim V^{g_2,g_3}\\
+\dim V^{g_2g_3}-\dim V^{g_2,g_3}-\dim V^{g_1,g_2g_3}+\dim V^{g_1,g_2,g_3}
\end{multline*}
respectively.  But these are easily seen to be equal, and the isomorphism \eqref{AssociativityIsomorphism} now follows from the discussion in Section~\ref{AgeObstructionSection}.
\end{proof}

\begin{theorem}[Chen-Hu \cite{\ChenHu}]\label{ChenHuTheorem} Let $X$ be an almost-complex abelian orbifold.  Then $E_{(g_1,g_2)}\to X_{(g_1,g_2)}$ is the summand of $\varepsilon^\ast TX$ spanned by the $\langle g_1,g_2\rangle$-invariant lines $L$ on which $\iota_L(g_1)+\iota_L(g_2)>1$.
\end{theorem}

\begin{proof}[Sketch Proof.]
By Theorem~\ref{BundleTheorem} it suffices to show that the functor $V\mapsto (V\otimes H^{0,1}_{\bar\partial}(\Sigma))^{\langle g_1,g_2 \rangle}$ is isomorphic to the functor that sends $V$ to the span of those linear subrepresentations $L\leqslant V$ for which $\iota_L(g_1)+\iota_L(g_2)>1$.  These functors are additive and so to prove this we must show that the corresponding additive functions obtained by computing the dimensions are equal, and to do this it suffices to check the claim for linear representations.  But it is trivial to verify that if $L$ is a linear representation of a finite group $G$ and $g_1,g_2\in G$, then
\[\iota_L(g_1)+\iota_L(g_2)-\iota_L(g_1g_2)-\dim L^{g_1g_2}+\dim L^{g_1,g_2}\]
is equal to $0$ if $\iota_L(g_1)+\iota_L(g_2)\leqslant 1$, and is equal to $1$ otherwise.  This completes the proof.
\end{proof}

\begin{theorem}[Gonz\'alez et al.~\cite{\GLSUX}]
For an almost-complex orbifold $X$ we have a ring-isomorphism $H^\ast_\mathrm{CR}(T^\ast X)\cong H^\ast_\mathrm{virt}(\Lambda X)$, where $H^\ast_\mathrm{virt}(\Lambda X)$ is the `virtual cohomology' of $\Lambda X$.
\end{theorem}

\begin{proof}[Sketch Proof.]
In \cite{\GLSUX} the proof was reduced using a cohomological argument to the claim that over a component $(T^\ast X)_{(g_1,g_2)}$ of $\Lambda^2(T^\ast X)$ we have 
\begin{equation}\label{CotangentEquation}E\oplus \pi^\ast\frac{\epsilon^\ast TX^{g_1g_2}}{\epsilon^\ast TX^{g_1,g_2}}\cong\pi^\ast\frac{\epsilon^\ast TX}{\epsilon^\ast TX^{g_1}+\epsilon^\ast TX^{g_2}},\end{equation}
where $\pi$ denotes the projection $T^\ast X\to X$ and the map it induces on $2$-sectors.

Each side of \eqref{CotangentEquation} is obtained by applying an additive functor to $\varepsilon^\ast TX$, so to prove \eqref{CotangentEquation} it suffices to show that the dimension of these functors when applied to a representation $V$ are always equal.  But the dimensions are
\begin{multline*}\iota_{V\oplus\bar V}(g_1)+\iota_{V\oplus\bar V}(g_2)-\iota_{V\oplus\bar V}(g_1g_2) -\dim(V\oplus\bar V)^{g_1g_2}+\dim(V\oplus\bar V)^{g_1,g_2} \\ +\dim V^{g_1g_2}-\dim V^{g_1,g_2}\end{multline*}
and
\[\dim V -\dim V^{g_1}-\dim V^{g_2}.\]
We have used the fact that $T(T^\ast X)\cong \pi^\ast TX\oplus\pi^\ast\overline{TX}$.  By noting that $\iota_{V\oplus\bar V}(g)=\iota_V(g)+\iota_{\bar V}(g)=\dim V-\dim V^g$ one verifies that the dimensions are equal, and this proves the theorem.
\end{proof}

\bibliographystyle{alpha}
\bibliography{AgeBibliography}

\newcommand{\etalchar}[1]{$^{#1}$}
\begin{thebibliography}{GLS{\etalchar{+}}07}

\bibitem[ALR07]{MR2359514}
Alejandro Adem, Johann Leida, and Yongbin Ruan.
\newblock {\em Orbifolds and stringy topology}, volume 171 of {\em Cambridge
  Tracts in Mathematics}.
\newblock Cambridge University Press, Cambridge, 2007.

\bibitem[ARZ06]{ARZ}
Alejandro Adem, Yongbin Ruan, and Bin Zhang.
\newblock Title: A stringy product on twisted orbifold k-theory, 2006.
\newblock Preprint, available at {\tt arXiv:0605534}.

\bibitem[CH06]{MR2242619}
Bohui Chen and Shengda Hu.
\newblock A de{R}ham model for {C}hen-{R}uan cohomology ring of abelian
  orbifolds.
\newblock {\em Math. Ann.}, 336(1):51--71, 2006.

\bibitem[CR04]{MR2104605}
Weimin Chen and Yongbin Ruan.
\newblock A new cohomology theory of orbifold.
\newblock {\em Comm. Math. Phys.}, 248(1):1--31, 2004.

\bibitem[FG03]{MR1971293}
Barbara Fantechi and Lothar G{\"o}ttsche.
\newblock Orbifold cohomology for global quotients.
\newblock {\em Duke Math. J.}, 117(2):197--227, 2003.

\bibitem[GLS{\etalchar{+}}07]{MR2318652}
Ana Gonz{\'a}lez, Ernesto Lupercio, Carlos Segovia, Bernardo Uribe, and
  Miguel~A. Xicot{\'e}ncatl.
\newblock Chen-{R}uan cohomology of cotangent orbifolds and {C}has-{S}ullivan
  string topology.
\newblock {\em Math. Res. Lett.}, 14(3):491--501, 2007.

\bibitem[IR96]{MR1463181}
Yukari Ito and Miles Reid.
\newblock The {M}c{K}ay correspondence for finite subgroups of {${\rm
  SL}(3,\bold C)$}.
\newblock In {\em Higher-dimensional complex varieties ({T}rento, 1994)}, pages
  221--240. de Gruyter, Berlin, 1996.

\bibitem[JKK07]{MR2285746}
Tyler~J. Jarvis, Ralph Kaufmann, and Takashi Kimura.
\newblock Stringy {$K$}-theory and the {C}hern character.
\newblock {\em Invent. Math.}, 168(1):23--81, 2007.

\end{thebibliography}
\end{document}